\documentclass[11pt, a4paper, twoside]{article}
\usepackage{amsfonts}
\usepackage{mathrsfs}
\usepackage{amsmath, amsthm}
\usepackage{pgf,tikz}
\usepackage{amssymb}
\usepackage{fancyhdr}
\usepackage[colorlinks, linkcolor=black, anchorcolor=black, citecolor=black]{hyperref}

\numberwithin{equation}{section}

\allowdisplaybreaks

\begin{document}

\fancyhf{}

\fancyhead[OR]{\thepage}

\renewcommand{\headrulewidth}{0pt}
\renewcommand{\thefootnote}{\fnsymbol {footnote}}

\theoremstyle{plain} 
\newtheorem{thm}{\indent\sc Theorem}[section] 
\newtheorem{lem}[thm]{\indent\sc Lemma}
\newtheorem{cor}[thm]{\indent\sc Corollary}
\newtheorem{prop}[thm]{\indent\sc Proposition}
\newtheorem{claim}[thm]{\indent\sc Claim}
\theoremstyle{definition} 
\newtheorem{dfn}[thm]{\indent\sc Definition}
\newtheorem{rem}[thm]{\indent\sc Remark}
\newtheorem{ex}[thm]{\indent\sc Example}
\newtheorem{notation}[thm]{\indent\sc Notation}
\newtheorem{assertion}[thm]{\indent\sc Assertion}
%
%
\numberwithin{equation}{section}
\renewcommand{\proofname}{\indent\sc Proof.} 
\def\C{\mathbb{C}}
\def\R{\mathbb{R}}
\def\Rn{{\mathbb{R}^n}}
\def\M{\mathbb{M}}
\def\N{\mathbb{N}}
\def\Q{{\mathbb{Q}}}
\def\Z{\mathbb{Z}}
\def\F{\mathcal{F}}
\def\L{\mathcal{L}}
\def\S{\mathcal{S}}
\def\supp{\operatorname{supp}}
\def\essi{\operatornamewithlimits{ess\,inf}}
\def\esss{\operatornamewithlimits{ess\,sup}}
\def\dlim{\displaystyle\lim}

\fancyhf{}

\fancyhead[EC]{W. LI, H. Wang, D. Yan}

\fancyhead[EL]{\thepage}

\fancyhead[OC]{Pointwise Convergence of Schr\"{o}dinger Means}

\fancyhead[OR]{\thepage}

\renewcommand{\headrulewidth}{0pt}
\renewcommand{\thefootnote}{\fnsymbol {footnote}}

\title{\textbf{Pointwise Convergence for sequences of Schr\"{o}dinger means in $\mathbb{R}^{2}$}
\footnotetext {This work is supported by the National Natural Science Foundation of China (No.11871452); Natural Natural Science Foundation of China (No.11701452); China Postdoctoral Science Foundation (No.2017M613193);  Natural Science Basic Research Plan in Shaanxi Province of China (No.2017JQ1009).}
\footnotetext {{}{2000 \emph{Mathematics Subject
 Classification}: 42B20, 42B25, 35S10.}}
\footnotetext {{}\emph{Key words and phrases}: Schr\"{o}dinger mean, Pointwise convergence. } } \setcounter{footnote}{0}
\author{
Wenjuan Li, Huiju Wang, Dunyan Yan}

\date{}
\maketitle

\begin{abstract}
We consider pointwise convergence of Schr\"{o}dinger means $e^{it_{n}\Delta}f(x)$ for $f \in H^{s}(\mathbb{R}^{2})$ and decreasing sequences $\{t_{n}\}_{n=1}^{\infty}$ converging to zero.
The main theorem  improves the previous results of [Sj\"{o}lin, JFAA, 2018] and [Sj\"{o}lin-Str\"{o}mberg, JMAA, 2020]  in $\mathbb{R}^{2}$. This study is based on investigating properties of Schr\"{o}dinger type maximal functions related to hypersurfaces with vanishing Gaussian curvature.
\end{abstract}

\section{Introduction}
The solution to the Schr\"{o}dinger equation
\begin{equation}\label{Eq11}
i{u_t} - \Delta u = 0, (x,t) \in {\mathbb{R}^N} \times \mathbb{R}^{+},
\end{equation}
with initial datum $u\left( {x,0} \right) = f,$ is formally written as
\[{e^{it\Delta }}f\left( x \right): = \int_{{\mathbb{R}^N}} {{e^{i\left( {x \cdot \xi  + t{{\left| \xi  \right|}^2}} \right)}}\widehat{f}} \left( \xi  \right)d\xi .\]
The problem about finding optimal $s$ for which
\begin{equation}\label{Eq12}
\mathop {\lim }\limits_{t \to 0^{+}} {e^{it\Delta }}f\left( x \right) = f(x) \hspace{0.2cm} a.e.\hspace{0.2cm} x\in \mathbb{R}^N,
\end{equation}
whenever $f \in {H^s}\left( {{\mathbb{R}^N}} \right),$ was first considered by Carleson \cite{C}, and extensively studied by Sj\"{o}lin \cite{S} and Vega \cite{V}, who proved independently the convergence for $s > 1/2$ in all dimensions. Dahlberg-Kenig \cite{DK} showed that the convergence does not hold for $s < 1/4$ in any dimension. In 2016, Bourgain \cite{B} gave conterexample showing that convergence can fail if $s<\frac{N}{2(N+1)}$. Very recently, Du-Guth-Li \cite{DGL} and Du-Zhang \cite{DZ} obtained the sharp results by the polynomial partitioning and decoupling method.

One of the natural generalizations of the pointwise convergence problem is to ask a.e. convergence of the Schr\"{o}dinger means where the limit is taken over decreasing sequences  $\{t_{n}\}_{n=1}^{\infty}$ converging to zero. That is to investigate relationship between optimal $s$  and properties of $\{t_{n}\}_{n=1}^{\infty}$  such that for each  function $f \in H^{s}(\mathbb{R}^{N})$,
 \begin{equation}\label{Eq1.3}
 \lim_{n \rightarrow \infty}e^{it_{n}\Delta}f(x) = f(x) \hspace{0.2cm} a.e.\hspace{0.2cm} x\in \mathbb{R}^N.
 \end{equation}
This problem was first considered by Sj\"{o}lin \cite{S1} in general dimensions and later improved by Sj\"{o}lin-Str\"{o}mberg \cite{SS1}. More recently, Dimou-Seeger \cite{DS} obtained a sharp characterization of this problem in the one-dimensional case. But in higer dimensional case, the sharp characterization of this problem still remains open.

In order to characterize the convergence of  $\{t_{n}\}_{n=1}^{\infty}$,  the Lorentz space ${\ell}^{r,\infty}(\mathbb{N})$, $r>0$ is involved. The sequence $\{t_{n}\}_{n=1}^{\infty} \in {\ell}^{r,\infty}(\mathbb{N})$ if and only  if
\begin{equation}
\mathop{sup}_{b>0}b^{r}\sharp\biggl\{n:t_{n}>b\biggl\} < \infty.
\end{equation}
 Dimou-Seeger \cite{DS} proved that when $N= 1$, $s  \ge min\{ \frac{r}{2r+1}, \frac{1}{4} \}$ is sufficient for (\ref{Eq1.3}) to hold. For $N>1$, it follows from \cite{S1} that $s >min\{ r, \frac{N}{2(N+1)} \}$ is sufficient. This was later improved by \cite{SS1} where $s >min\{ \frac{r}{1+r}, \frac{N}{2(N+1)} \}$ is shown to be enough for pointwise convergence. The main theorem of this paper improves the previous  results of \cite{S1} and \cite{SS1} in dimension two.

By standard arguments, in order to obtain the convergence result, it is sufficient to
show the maximal function estimate in $\mathbb{R}^{2}$. Our main results are as follows.

\begin{thm}\label{theorem1.1}
Given a decreasing sequence $\{t_{n}\}_{n=1}^{\infty} \in {\ell}^{r,\infty}(\mathbb{N})$ converging to zero and $\{t_{n}\}_{n=1}^{\infty} \subset (0,1)$,
then for any $s> s_{0}=min\{\frac{r}{\frac{4}{3}r+1}, \frac{1}{3}\}$, we have
\begin{equation}\label{Eq1.7+}
\biggl\|\mathop{sup}_{n \in \mathbb{N}} |e^{it_{n}\Delta}f|\biggl\|_{L^{2}(B(0,1))} \leq C\|f\|_{H^s(\mathbb{R}^2)},
\end{equation}
whenever $f\in H^s(\mathbb{R}^2)$, where the constant $C$  does not depend on $f$.
\end{thm}

By translating invariance in the $x-$variables, $B(0,1)$ in Theorem \ref{theorem1.1} can be replaced by any ball of radius $1$ in $\mathbb{R}^{2}$. Then we obtain the following convergence result.

\begin{thm}\label{theorem1.2}
Given a decreasing sequence $\{t_{n}\}_{n=1}^{\infty} \in {\ell}^{r,\infty}(\mathbb{N})$ converging to zero and $\{t_{n}\}_{n=1}^{\infty} \subset (0,1)$,
then for any $s> s_{0}=min\{\frac{r}{\frac{4}{3}r+1}, \frac{1}{3}\}$,  (\ref{Eq1.3}) holds
whenever $f\in H^s(\mathbb{R}^2)$.
\end{thm}

Theorem \ref{theorem1.2} improves the previous results of \cite{S1} and \cite{SS1} when $N=2$ since $\frac{r}{\frac{4}{3}r+1} < \frac{r}{r+1}$ for $r >0$. We  pose two examples for $\{t_{n}\}_{n=1}^{\infty}$. It is not hard to check that  $\{t_{n}\}_{n=1}^{\infty} \in {\ell}^{r,\infty}(\mathbb{N})$  if we take (E1): $t_{n}=\frac{1}{n^{1/r}}$, $n \ge1$. It is obvious that when $r < \frac{3}{5}$, there is a gain over the general pointwise convergence result for $s>1/3$ in $\mathbb{R}^2$.
 Another example is the lacunary sequence (E2): $t_{n} = 2^{-n}$, $n \ge 1$. For this example, it is worth to mention that  $\{t_{n}\}_{n=1}^{\infty} \in {\ell}^{r,\infty}(\mathbb{N})$ for each $r >0$. Therefore, inequality (\ref{Eq1.3}) holds whenever $f\in H^s(\mathbb{R}^2)$ for any $s >0$.

We briefly sketch the proof of Theorem \ref{theorem1.1}, and leave the details to Section 2. Notice that when $\frac{r}{\frac{4}{3}r+1} \ge \frac{1}{3}$,  Theorem \ref{theorem1.1} follows from H\"{o}lder's inequality and the  theorem below.
\begin{thm}\label{lemma3.2}(\cite{DGL})
 For any $s>1/3$, the following bound holds: for any function $f\in H^s(\mathbb{R}^2)$,
\begin{equation*}
\biggl\|\mathop{sup}_{0<t<1}|e^{it\Delta}f(x)|\biggl\|_{L^3(B(0,1))} \leq C_s\|f\|_{H^s(\mathbb{R}^2)}.
\end{equation*}
\end{thm}
Therefore, we only need to consider the case when
$\frac{r}{\frac{4}{3}r+1} < \frac{1}{3}.$

By Littlewood-Paley decomposition, it is not hard to deal with the low frequency parts by standard argument. So next we just concentrate ourselves on the case when $\textmd{supp}\hat{f} \subset  \{\xi: |\xi| \sim 2^{k}\}$, $k \gg 1$. We decompose $\{t_{n}\}_{n=1}^{\infty}$ as
\[A^{1}_{k}:= \biggl\{t_{n}: t_{n} \ge 2^{-\frac{2k}{\frac{4}{3}r+1}} \biggl\}\]
and
\[A^{2}_{k}:= \biggl\{t_{n}: t_{n} < 2^{-\frac{2k}{\frac{4}{3}r+1}} \biggl\}.\]
Then we consider the maximal function
\[\mathop{sup}_{n \in \mathbb{N}: t_{n} \in A^{1}_{k}} |e^{it_{n}\Delta}f|\]
and
\[\mathop{sup}_{n \in \mathbb{N}: t_{n} \in A^{2}_{k}} |e^{it_{n}\Delta}f|,\]
respectively. We deal with the first term by the assumption that the decreasing sequence $\{t_{n}\}_{n=1}^{\infty} \in {\ell}^{r,\infty}(\mathbb{N})$ and Plancherel theorem. For the second term, since we only need to consider the case when $\frac{r}{\frac{4}{3}r+1} < \frac{1}{3}$, we may assume that $0<r < \frac{3}{5}$, which yields $A_{k}^{2} \subset (0,2^{-\frac{2k}{\frac{4}{3}r+1}})$ for $k< \frac{2k}{\frac{4}{3}r+1}< 2k$. Hence we  reduce our proof to the following theorem.

\begin{thm}\label{theorem1.4}
If supp $\hat{f} \subset  \{\xi: |\xi| \sim 2^{k}\}$, then for any small constant $\epsilon >0$ and interval $I = (0,2^{-j})$, where
\[k < j < 2k, \:\ j \in \mathbb{R},\]
we have
 \begin{equation}\label{Eq1.8}
\biggl\|\mathop{sup}_{t \in I}|e^{it\Delta}f(x)|\biggl\|_{L^2(B(0,1))} \leq C2^{(k-\frac{j}{2})\frac{3}{4} +\epsilon k}\|f\|_{L^{2}(\mathbb{R}^2)},
\end{equation}
where the constant $C$ does not depend on $f$.
\end{thm}

We will prove Theorem \ref{theorem1.4} in Section 3. In one-dimensional case, similar result was built in \cite{DS} by $TT^{\star}$ and stationary phase method. But their method seems not to work well in the higher dimensional case. If $j \ge 2k$ in Theorem \ref{theorem1.4}, then the length of $I$ is small enough so that we can apply Sobolev's embedding to obtain
 \begin{equation}\label{uniform}
\biggl\|\mathop{sup}_{t \in I}|e^{it\Delta}f(x)|\biggl\|_{L^2(B(0,1))} \leq C\|f\|_{L^{2}(\mathbb{R}^2)}.
\end{equation}
When $j=k$, Lee \cite{L} applied the bilinear method to show that
 \begin{equation}
\biggl\|\mathop{sup}_{t \in I}|e^{it\Delta}f(x)|\biggl\|_{L^2(B(0,1))} \leq C2^{\frac{3}{8}k + \epsilon k}\|f\|_{L^{2}(\mathbb{R}^2)}.
\end{equation}
In our case, the exponent $2^{(k-\frac{j}{2})\frac{3}{4} +\epsilon k}$ is in the middle of $1$ and $2^{\frac{3}{8}k + \epsilon k}$ when $k<j<2k$.
However, inequality (1.8) was improved by \cite{DGL} and \cite{DZ} using Broad-Narrow argument and polynomial partitioning. It was proved that if $j=k$, then
 \begin{equation}\label{Eq1.9}
\biggl\|\mathop{sup}_{t \in I}|e^{it\Delta}f(x)|\biggl\|_{L^2(B(0,1))} \leq C2^{\frac{1}{3}k + \epsilon k}\|f\|_{L^{2}(\mathbb{R}^2)}.
\end{equation}
In view of technical difficulties, we could not improve our result to this exponent. Moreover,  notice that (\ref{Eq1.9}) also holds for $j<k$ due to the localizing lemma in Lee-Rogers \cite{LR}.

In our case $k<j<2k$, by scaling, we actually have to treat the maximal function defined by
\[\mathop{sup}_{t \in (0,1)}\biggl|\int_{\mathbb{R}^{2}}e^{i2^{k}(x\cdot \eta +2^{k-j}t|\eta|^{2})}\hat{g}(\eta)d\eta\biggl|,\]
where $\hat{g}$ is supported in the annular $\{\eta: |\eta| \sim 1\}$. It is obvious that the Gaussian curvature of the hypersurface $(\eta, 2^{k-j}|\eta|^{2})$ vanishes as $k$ tends to infinity. Therefore, it seems difficult to apply the classical analysis of local smoothing estimate or Broad-Narrow argument directly. Inspired by \cite{LW}, in which the authors established $L^{p}\rightarrow L^{q}$ estimates for Fourier integral operators related to hypersurfaces with vanishing Gaussian curvature, we try to use Whitney type decomposition and a bilinear result from \cite{L} to prove Theorem \ref{theorem1.4}. Theorem \ref{theorem1.5} below follows from Proposition 3.1 in \cite{L} and rescaling.

\begin{thm}\label{theorem1.5}(\cite{L})
Let $\lambda \gg 1$,  $l$: $2^{l} \ll \lambda^{\frac{1}{4}}$. If  $\hat{g_{1}}$, $\hat{g_{2}}$ are supported in $B(\eta_{0}, 2^{-l}) \subset \{\eta: |\eta| \sim 1\}$ and $dist(\text{supp }\widehat{g_{1}}, \text{ supp } \widehat{g_{2}}) \sim 2^{-l}$, then for any $\epsilon >0$, we have
\begin{equation}\label{Eqglobal}
\biggl\| \mathop{sup}_{t \in (0,1)}|T_{\lambda}g_{1} T_{\lambda}g_{2}|\biggl\|_{_{L^{2}(B(0,1))}} \leq C 2^{-\frac{l}{2}} \lambda^{\frac{3}{4}+\epsilon} \lambda^{-2}\|g_{1}\|_{L^2}\|g_{2}\|_{L^2},
\end{equation}
where for $i=1,2$,
\[T_{\lambda}g_{i}(x,t)=\int_{\mathbb{R}^{2}}e^{i\lambda(x\cdot \eta + t|\eta|^{2})}\hat{g_{i}}(\eta)d\eta.\]
\end{thm}

\textbf{Conventions}: Throughout this article, we shall use the well known notation $A\gg B$, which means if there is a sufficiently large constant $G$, which does not depend on the relevant parameters arising in the context in which
the quantities $A$ and $B$ appear, such that $ A\geq GB$. We write $A\sim B$, and mean that $A$ and $B$ are comparable. By
$A\lesssim B$ we mean that $A \le CB $ for some constant $C$ independent of the parameters related to  $A$ and $B$.

\section{Proof of Theorem \ref{theorem1.1}}
\textbf{Proof of Theorem \ref{theorem1.1}.}
Set
\[s_{1}=\frac{r}{\frac{4}{3}r+1} + \epsilon\]
for some sufficiently small constant $\epsilon >0$.
We decompose $f$ as
\[f=\sum_{k=0}^{\infty}{f_{k}},\]
where $\textmd{supp} \hat{f_{0}} \subset B(0,1)$, $\textmd{supp} \hat{f_{k}} \subset \{\xi: |\xi| \sim 2^{k}\}, k \ge 1$. Then we have
\begin{equation}\label{Eq2.2}
\biggl\|\mathop{sup}_{n \in \mathbb{N}} |e^{it_{n}\Delta}f|\biggl\|_{L^{2}(B(0,1))} \le \sum_{k=0}^{\infty}{\biggl\|\mathop{sup}_{n \in \mathbb{N}} |e^{it_{n}\Delta}f_{k}|\biggl\|_{L^{2}(B(0,1))}}.
\end{equation}

For $k \lesssim 1$, since for each $x \in B(0,1)$,
\[\biggl|e^{it_{n}\Delta}f_{k}(x)\biggl| \lesssim \|f_{k}\|_{L^{2}(\mathbb{R}^{2})},\]
it is obvious that
\begin{equation}\label{Eq2.3}
\biggl\|\mathop{sup}_{n \in \mathbb{N}} |e^{it_{n}\Delta}f_{k}|\biggl\|_{L^{2}(B(0,1))} \lesssim \|f\|_{H^{s_{1}}(\mathbb{R}^{2})}.
\end{equation}

For each $k \gg 1$, we decompose $\{t_{n}\}_{n=1}^{\infty}$ as
\[A^{1}_{k}:= \biggl\{t_{n}: t_{n} \ge 2^{-\frac{2k}{\frac{4}{3}r+1}} \biggl\}\]
and
\[A^{2}_{k}:= \biggl\{t_{n}: t_{n} < 2^{-\frac{2k}{\frac{4}{3}r+1}} \biggl\}.\]
Then we have
\begin{align}\label{Eq2.4}
\biggl\|\mathop{sup}_{n \in \mathbb{N}} |e^{it_{n}\Delta}f_{k}|\biggl\|_{L^{2}(B(0,1))} &\le \biggl\|\mathop{sup}_{n \in \mathbb{N}: t_{n} \in A^{1}_{k}} |e^{it_{n}\Delta}f_{k}|\biggl\|_{L^{2}(B(0,1))} \nonumber\\
& \:\ + \biggl\|\mathop{sup}_{n \in \mathbb{N}: t_{n} \in A^{2}_{k}} |e^{it_{n}\Delta}f_{k}|\biggl\|_{L^{2}(B(0,1))} \nonumber\\
&: =I + II.
\end{align}

We firstly estimate $I$.
Since  $\{t_{n}\}_{n=1}^{\infty} \in {\ell}^{r,\infty}(\mathbb{N})$ and assumption $r\in (0,3/5)$, we have
\begin{equation}\label{Eq2.5}
\sharp A^{1}_{k}  \le C 2^{\frac{2rk}{\frac{4}{3}r+1}},
\end{equation}
which implies that
\begin{align}\label{Eq2.6}
I &\le \biggl(\sum_{n \in \mathbb{N}: t_{n} \in A^{1}_{k}}{ \biggl\|e^{it_{n}\Delta}f_{k}\biggl\|^{2}_{L^{2}(B(0,1))}}\biggl)^{1/2} \nonumber\\
&\le 2^{\frac{rk}{\frac{4}{3}r+1}} \|f_{k}\|_{L^{2}(\mathbb{R}^{2})} \nonumber\\
&\lesssim 2^{-\epsilon k}\|f\|_{H^{s_1}(\mathbb{R}^{2})}.
\end{align}

For $II$, since
\[A^{2}_{k} \subset \biggl(0,  2^{-\frac{2k}{\frac{4}{3}r+1}}\biggl).\]
 By previous discussion, we have
$k<\frac{2k}{\frac{4}{3}r+1}<2k.$
Then it follows from Theorem \ref{theorem1.4} that,
\begin{equation}\label{Eq2.7}
II \lesssim 2^{(\frac{r}{\frac{4}{3}r+1}+\frac{\epsilon}{2})k}\|f_{k}\|_{L^{2}(\mathbb{R}^{2})} \le 2^{-\frac{\epsilon}{2} k}\|f\|_{H^{s_1}(\mathbb{R}^{2})}.
\end{equation}

Inequalities (\ref{Eq2.4}), (\ref{Eq2.6}) and (\ref{Eq2.7}) yield for $k \gg 1$,
\begin{align}\label{Eq2.8}
\biggl\|\mathop{sup}_{n \in \mathbb{N}} |e^{it_{n}\Delta}f_{k}|\biggl\|_{L^{2}(B(0,1))} &\lesssim 2^{-\frac{\epsilon k}{2}}\|f\|_{H^{s_1}(\mathbb{R}^{2})}.
\end{align}

Combining  (\ref{Eq2.2}), (\ref{Eq2.3}) and (\ref{Eq2.8}),  inequality (\ref{Eq1.7+}) holds true for $s_1$. By the arbitrariness of $\epsilon$, in fact, we can get for any $s>s_0$, inequality (\ref{Eq1.7+}) remains true.


\section{Proof of Theorem \ref{theorem1.4}}
\textbf{Proof of Theorem \ref{theorem1.4}.} If $2k-k\epsilon\leq j<2k$, then $I\subset (0,2^{-2k+k\epsilon})$, which can be covered by $2^{k\epsilon}$ intervals of length $2^{-2k}$, from inequality (\ref{uniform}), we have that
\begin{equation}
\biggl\|\mathop{sup}_{t \in I}|e^{it\Delta}f(x)|\biggl\|_{L^2(B(0,1))}\leq 2^{k\epsilon}\|f\|_{L^2}.
\end{equation}
Next, we always assume that $k<j<2k-k\epsilon$. By Whitney type decomposition, we have
\[(e^{it\Delta}f)^{2}= \sum_{l:  1  \le 2^{l} \le 2^{(2k-j)/4 }} \sum_{m,m^{\prime}: dist(\text{supp }\widehat{f_{m}^{l}}, \text{ supp } \widehat{ f_{m^{\prime}}^{l}}) \sim 2^{k-l}}  e^{it\Delta}f_{m}^{l} \cdot e^{it\Delta}f_{m^{\prime}}^{l},\]
where $l \in \mathbb{N}$, $\text{supp }\widehat{f_{m}^{l}}$ and $ \text{supp } \widehat{ f_{m^{\prime}}^{l}}$ are contained in cubes with side length $2^{k-l}$. When $2^{l}  \sim 2^{\frac{2k-j}{4}}$, we abuse the notation that by saying $dist(\text{supp }\widehat{f_{m}^{l}}, \text{ supp } \widehat{ f_{m^{\prime}}^{l}}) \sim 2^{k-l}$ we mean $dist(\text{supp }\widehat{f_{m}^{l}}, \text{supp } \widehat{ f_{m^{\prime}}^{l}}) \lesssim  2^{k-l}$. Then we have
 \begin{align}\label{Eq3.2}
&\biggl\|\mathop{sup}_{t \in I}|e^{it\Delta}f(x)|\biggl\|_{L^2(B(0,1))} \nonumber\\
&\le \biggl\{\sum_{l:  1  \le 2^{l} \le 2^{(2k-j)/4 } } \sum_{m,m^{\prime}: dist(\text{supp }\widehat{f_{m}^{l}}, \text{ supp } \widehat{ f_{m^{\prime}}^{l}}) \sim 2^{k-l}} \biggl\|\mathop{sup}_{t \in I}|e^{it\Delta}f_{m}^{l}| | e^{it\Delta}f_{m^{\prime}}^{l}| \biggl\|_{L^1(B(0,1))} \biggl \}^{1/2}.
\end{align}

It is sufficient to prove the following two lemmas.

\begin{lem}\label{main}
For each $l$: $2^{l} \ll 2^{(2k-j)/4}$, we have
\begin{equation*}
\biggl\|\mathop{sup}_{t \in I}|e^{it\Delta}f_{m}^{l}| | e^{it\Delta}f_{m^{\prime}}^{l}| \biggl\|_{L^1(B(0,1))} \leq C 2^{-\frac{l}{2}}2^{(2k-j)\frac{3}{4} +2\epsilon k} \|f_{m}^{l}\|_{L^2}\|f_{m^{\prime}}^{l}\|_{L^2}.
\end{equation*}
\end{lem}

\begin{lem}\label{basic}
For each $l$: $2^{l} \sim 2^{(2k-j)/4}$, we have
\begin{equation*}
\biggl\|\mathop{sup}_{t \in I}|e^{it\Delta}f_{m}^{l}| | e^{it\Delta}f_{m^{\prime}}^{l}| \biggl\|_{L^1(B(0,1))} \leq C 2^{(2k-j)\frac{3}{4}} \|f_{m}^{l}\|_{L^2}\|f_{m^{\prime}}^{l}\|_{L^2}.
\end{equation*}
\end{lem}

Indeed, if Lemma \ref{main} and Lemma \ref{basic} hold true, then by the Cauchy-Schwartz inequality, we obtain
 \begin{align}\label{Eq3.2}
&\biggl\|\mathop{sup}_{t \in I}|e^{it\Delta}f(x)|\biggl\|_{L^2(B(0,1))} \nonumber\\
&\le \biggl\{\sum_{l:   1  \le 2^{l} \ll 2^{(2k-j)/4 } }C 2^{-\frac{l}{2}}2^{(2k-j)\frac{3}{4} +2\epsilon k} \sum_{m,m^{\prime}: dist(\text{supp }\widehat{f_{m}^{l}}, \text{ supp } \widehat{ f_{m^{\prime}}^{l}}) \sim 2^{k-l}} \|f_{m}^{l}\|_{L^2}\|f_{m^{\prime}}^{l}\|_{L^2} \biggl \}^{1/2} \nonumber\\
& \:\ + \biggl\{\sum_{l:  2^{l} \sim 2^{(2k-j)/4} }C 2^{(2k-j)\frac{3}{4} } \sum_{m,m^{\prime}: dist(\text{supp }\widehat{f_{m}^{l}}, \text{ supp } \widehat{ f_{m^{\prime}}^{l}}) \sim 2^{k-l}} \|f_{m}^{l}\|_{L^2}\|f_{m^{\prime}}^{l}\|_{L^2} \biggl \}^{1/2} \nonumber\\
&\lesssim 2^{(k-\frac{j}{2})\frac{3}{4} +\epsilon k}\|f\|_{L^{2}},
\end{align}
which yields Theorem \ref{theorem1.4}.

Let's turn to prove Lemma \ref{main} and Lemma \ref{basic}. Lemma \ref{basic} follows from H\"{o}lder's inequality and the following two estimates
\begin{align}\label{Eq3.3}
\biggl\|\mathop{sup}_{t \in I}|e^{it\Delta}f_{m}^{l}| \biggl\|_{L^2(B(0,1))} \leq C 2^{(k-\frac{j}{2})\frac{3}{4}} \|f_{m}^{l}\|_{L^2},
\end{align}
\begin{align}\label{Eq3.4}
\biggl\|\mathop{sup}_{t \in I}|e^{it\Delta}f_{m^{\prime}}^{l}| \biggl\|_{L^2(B(0,1))} \leq C 2^{(k-\frac{j}{2})\frac{3}{4}} \|f_{m^{\prime}}^{l}\|_{L^2}.
\end{align}
We only prove inequality (\ref{Eq3.3}) since the proof of inequality (\ref{Eq3.4}) is similar. Without loss of generality, we may assume that $\text{supp }\widehat{f_{m}^{l}} \subset B(\xi_{0}, 2^{\frac{2k+j}{4}})$, where $|\xi_{0}| \lesssim 2^{k}$.
By changes of variables,
\[\xi = \zeta + \xi_{0}, \:\ |\zeta| \le 2^{\frac{2k+j}{4}},\]
 we have
 \[\mathop{sup}_{t \in I}|e^{it\Delta}f_{m}^{l}| = \mathop{sup}_{t \in I}\biggl| \int_{{\mathbb{R}^2}} {{e^{i[ {x \cdot \zeta  + t({{\left| \zeta  \right|}^2}} +2 \zeta \cdot \xi_{0}) ]}}\widehat{f_{m}^{l}}} \left( \zeta + \xi_{0}  \right)d\zeta \biggl| .\]
 It follows from Sobolev's embedding and Plancherel theorem that
 \begin{align}
 &\biggl\|\mathop{sup}_{t \in I}|e^{it\Delta}f_{m}^{l}| \biggl\|_{L^2(B(0,1))} \nonumber\\
 &\le \|f_{m}^{l}\|_{L^2} +  \biggl\|  \int_{{\mathbb{R}^2}} {{e^{i[ {x \cdot \zeta  + t({{\left| \zeta  \right|}^2}} +2 \zeta \cdot \xi_{0}) ]}}\widehat{f_{m}^{l}}} \left( \zeta + \xi_{0}  \right)d\zeta  \biggl\|^{1/2}_{L^2(B(0,1)\times I)} \nonumber\\
 & \:\  \times  \biggl\|  \int_{{\mathbb{R}^2}} {{e^{i[ {x \cdot \zeta  + t({{\left| \zeta  \right|}^2}} +2 \zeta \cdot \xi_{0}) ]}} (|\zeta|^{2} + 2\zeta \cdot \xi_{0})\widehat{f_{m}^{l}}} \left( \zeta + \xi_{0}  \right)d\zeta  \biggl\|^{1/2}_{L^2(B(0,1)\times I)} \nonumber\\
 &\le \|f_{m}^{l}\|_{L^2} + 2^{-\frac{j}{2}} \biggl\|\widehat{f_{m}^{l}}(\zeta + \xi_{0})\biggl\|_{L^{2}}^{1/2}\biggl\| (|\zeta|^{2} + 2\zeta \cdot \xi_{0})\widehat{f_{m}^{l}}(\zeta + \xi_{0})\biggl\|_{L^{2}}^{1/2} \nonumber\\
 &\le \|f_{m}^{l}\|_{L^2} + 2^{-\frac{j}{2}}2^{\frac{k}{2}} 2^{\frac{2k+j}{8}}\biggl\|\widehat{f_{m}^{l}}(\zeta + \xi_{0})\biggl\|_{L^{2}} \nonumber\\
 &\leq C 2^{(k-\frac{j}{2})\frac{3}{4}} \|f_{m}^{l}\|_{L^2}. \nonumber
 \end{align}
Then we arrive at inequality (\ref{Eq3.3}).

We will prove Lemma \ref{main} in the rest of this  section. By rescaling, we turn to estimate
\begin{align}\label{estimate:L4+}
&\mathop{sup}_{t \in I}|e^{it\Delta}f_{m}^{l}(x)| | e^{it\Delta}f_{m^{\prime}}^{l}(x)| \nonumber\\
&=\mathop{sup}_{t \in I} \biggl| \int_{{\mathbb{R}}^2}e^{i2^{k}(\eta \cdot x + 2^{k}t|\eta|^{2})}2^{2k}\widehat{f_{m}^{l}}(2^{k}\eta)d\eta \biggl| \biggl| \int_{{\mathbb{R}}^2}e^{i2^{k}(\eta^{\prime} \cdot x + 2^{k}t|\eta^{\prime}|^{2})}2^{2k}\widehat{f_{m^{\prime}}^{l}}(2^{k}\eta^{\prime})d\eta^{\prime} \biggl| \nonumber\\
&=\mathop{sup}_{t \in (0,1)} \biggl| \int_{{\mathbb{R}}^2}e^{i2^{k}(\eta \cdot x + 2^{k-j}t|\eta|^{2})}2^{2k}\widehat{f_{m}^{l}}(2^{k}\eta)d\eta \biggl| \biggl| \int_{{\mathbb{R}}^2}e^{i2^{k}(\eta ^{\prime}\cdot x + 2^{k-j}t|\eta^{\prime}|^{2})}2^{2k}\widehat{f_{m^{\prime}}^{l}}(2^{k}\eta^{\prime})d\eta^{\prime} \biggl| \nonumber\\
&=\mathop{sup}_{t \in (0,1)} \biggl| \int_{{\mathbb{R}}^2}e^{i2^{k}(\eta \cdot x + 2^{k-j}t|\eta|^{2})}\widehat{F_{m}^{l}}(\eta)d\eta \biggl| \biggl| \int_{{\mathbb{R}}^2}e^{i2^{k}(\eta^{\prime} \cdot x + 2^{k-j}t|\eta^{\prime}|^{2})}\widehat{F_{m^{\prime}}^{l}}(\eta^{\prime})d\eta^{\prime} \biggl|,
\end{align}
where
\[\widehat{F_{m}^{l}}(\eta)= 2^{2k}\widehat{f_{m}^{l}}(2^{k}\eta),\]
\[\widehat{F_{m^{\prime}}^{l}}(\eta^{\prime})= 2^{2k}\widehat{f_{m^{\prime}}^{l}}(2^{k}\eta^{\prime}).\]
Here we notice that $\text{supp }\widehat{F_{m}^{l}}, \text{ supp } \widehat{ F_{m^{\prime}}^{l}}$ are contained in $\{\eta :|\eta| \sim 1\}$. More concretely, $\text{supp }\widehat{F_{m}^{l}}, \text{ supp } \widehat{ F_{m^{\prime}}^{l}}$ are contained in cubes with side length $2^{-l}$ and $dist(\text{supp }\widehat{F_{m}^{l}}, \text{ supp } \widehat{ F_{m^{\prime}}^{l}}) \thicksim 2^{-l}$.

Next we will try to localize $x$ into cubes with side length $2^{k-j}$. We have
\begin{align}
\int_{{\mathbb{R}}^2}e^{i2^{k}(\eta \cdot x + 2^{k-j}t|\eta|^{2})}\widehat{F_{m}^{l}}(\eta)d\eta = \int_{\mathbb{R}^{2}} \int_{{\mathbb{R}}^2}e^{i2^{k}(\eta \cdot x + 2^{k-j}t|\eta|^{2})-iz \cdot \eta} \phi(\eta) d\eta F_{m}^{l}(z)dz, \nonumber
\end{align}
where $\phi \in C_{c}^{\infty}(\mathbb{R}^2)$ such that $\phi(\eta) = 1$ on $\{\eta: |\eta| \sim 1\}$ and decays rapidly outside.
Denote
\begin{align}
K(x,z,t)=  \int_{{\mathbb{R}}^2}e^{i2^{k}(\eta \cdot x + 2^{k-j}t|\eta|^{2})-iz \cdot \eta} \phi(\eta) d\eta . \nonumber
\end{align}
Integrating by parts shows that if $|x-z/2^{k}| \gg 2^{k-j}$,  then for each $t \in (0,1)$ and any positive integer $\tilde{M}\gg1$,
\begin{align}
|K(x,z,t)| \leq  \frac {C_{\tilde{M}+2}}{\biggl(1+ 2^{k} \biggl| |x-\frac{z}{2^{k}}| +\mathcal{O}(2^{k-j}) \biggl|   \biggl)^{\tilde{M}+2}}. \nonumber
\end{align}
Notice that we may assume  $2^{j} \le 2^{2k-k\epsilon}$ since  (\ref{uniform}) holds when $j=2k$, then there is a  sufficiently large $M<\tilde{M}$ such that
\begin{align}
 \int_{z:|x-z/2^{k}| \gg 2^{k-j}} K(x,z,t) F_{m}^{l}(z)dz \le C2^{-kM}\|F_{m}^{l}\|_{L^{2}}.
\end{align}
 By the same argument, if $|x-z^{\prime}/2^{k}| \gg 2^{k-j}$,  then
 \begin{align}
 \int_{z^{\prime}:|x-z^{\prime} /2^{k}| \gg 2^{k-j}} K(x,z^{\prime},t) F_{m^{\prime}}^{l}(z^{\prime})dz^{\prime} \le C2^{-kM}\|F_{m^{\prime}}^{l}\|_{L^{2}},
\end{align}
where
\begin{align}
K(x,z^{\prime},t)=  \int_{{\mathbb{R}}^2}e^{i2^{k}(\eta^{\prime} \cdot x + 2^{k-j}t|\eta^{\prime}|^{2})-iz^{\prime} \cdot \eta^{\prime}} \phi(\eta^{\prime}) d\eta^{\prime}. \nonumber
\end{align}

Let $\pi$ be a smooth function such that $\sum_{h\in \mathbb{Z}^2}\pi(\cdot-h)=1$ and its Fourier transform is supported in $B(0,1)$. Let $\{Q\}$ be a collection of cubes of side length $2^{k-j}$ which partitions  $\mathbb{R}^2$ and let $a_Q$ be the affine map sending $Q$ to the unit cube centered at the origin. Set $\Pi_Q=\pi\circ a_Q$, here $\Pi_Q$ is a smooth function essentially supported in $Q$.
Hence,
\begin{align}\label{sum}
&\biggl| \int_{{\mathbb{R}}^2}e^{i2^{k}(\eta \cdot x + 2^{k-j}t|\eta|^{2})}\widehat{F_{m}^{l}}(\eta)d\eta \biggl| \biggl| \int_{{\mathbb{R}}^2}e^{i2^{k}(\eta^{\prime} \cdot x + 2^{k-j}t|\eta^{\prime}|^{2})}\widehat{F_{m^{\prime}}^{l}}(\eta^{\prime})d\eta^{\prime} \biggl| \nonumber\\
&=\biggl| \int_{\mathbb{R}^{2}} K(x,z,t) \biggl( \sum_{Q}\Pi_{Q}(\frac{z}{2^{k}})  \biggl) F_{m}^{l}(z)dz \biggl|  \biggl| \int_{\mathbb{R}^{2}} K(x,z^{\prime},t)  \biggl( \sum_{Q^{\prime}} \Pi_{Q^{\prime}}  (\frac{z^{\prime}}{2^{k}}) \biggl) F_{m^{\prime}}^{l}(z^{\prime})dz^{\prime} \biggl| \nonumber\\
&=\biggl| \int_{\mathbb{R}^{2}} K(x,z,t) \biggl( \sum_{Q: dist(x,Q) \lesssim 2^{k-j}}\Pi_{Q}(\frac{z}{2^{k}})  \biggl) F_{m}^{l}(z)dz \biggl|   \nonumber\\
&\:\ \times \biggl| \int_{\mathbb{R}^{2}} K(x,z^{\prime},t) \biggl( \sum_{Q^{\prime}: dist (x,Q^{\prime}) \lesssim 2^{k-j}} \Pi_{Q^{\prime}}  (\frac{z^{\prime}}{2^{k}}) \biggl) F_{m^{\prime}}^{l}(z^{\prime})dz^{\prime} \biggl| \nonumber\\
&+\biggl| \int_{\mathbb{R}^{2}} K(x,z,t) \biggl( \sum_{Q: dist(x,Q) \lesssim 2^{k-j}}\Pi_{Q}(\frac{z}{2^{k}})  \biggl) F_{m}^{l}(z)dz \biggl|   \nonumber\\
&\:\ \times \biggl| \int_{\mathbb{R}^{2}} K(x,z^{\prime},t) \biggl( \sum_{Q^{\prime}: dist (x,Q^{\prime}) \gg 2^{k-j}} \Pi_{Q^{\prime}}  (\frac{z^{\prime}}{2^{k}}) \biggl) F_{m^{\prime}}^{l}(z^{\prime})dz^{\prime} \biggl| \nonumber\\
&+\biggl| \int_{\mathbb{R}^{2}} K(x,z,t) \biggl( \sum_{Q: dist(x,Q) \gg 2^{k-j}}\Pi_{Q}(\frac{z}{2^{k}})  \biggl) F_{m}^{l}(z)dz \biggl|   \nonumber\\
&\:\ \times \biggl| \int_{\mathbb{R}^{2}} K(x,z^{\prime},t) \biggl( \sum_{Q^{\prime}: dist (x,Q^{\prime}) \lesssim 2^{k-j}} \Pi_{Q^{\prime}}  (\frac{z^{\prime}}{2^{k}}) \biggl) F_{m^{\prime}}^{l}(z^{\prime})dz^{\prime} \biggl| \nonumber\\
&+\biggl| \int_{\mathbb{R}^{2}} K(x,z,t) \biggl( \sum_{Q: dist(x,Q) \gg 2^{k-j}}\Pi_{Q}(\frac{z}{2^{k}})  \biggl) F_{m}^{l}(z)dz \biggl|   \nonumber\\
&\:\ \times \biggl| \int_{\mathbb{R}^{2}} K(x,z^{\prime},t) \biggl( \sum_{Q^{\prime}: dist (x,Q^{\prime}) \gg 2^{k-j}} \Pi_{Q^{\prime}}  (\frac{z^{\prime}}{2^{k}}) \biggl) F_{m^{\prime}}^{l}(z^{\prime})dz^{\prime} \biggl| \nonumber\\
&\le \sum_{Q,Q^{\prime}: cQ \cap cQ^{\prime} \neq \emptyset} \Pi_{cQ \cap cQ^{\prime}}(x) \biggl| \int_{\mathbb{R}^{2}} K(x,z,t) \Pi_{Q}(\frac{z}{2^{k}}) F_{m}^{l}(z)dz \biggl|  \biggl| \int_{\mathbb{R}^{2}} K(x,z^{\prime},t) \Pi_{Q^{\prime}}(\frac{z^{\prime}}{2^{k}}) F_{m^{\prime}}^{l}(z^{\prime})dz^{\prime} \biggl| \nonumber\\
&\quad +C 2^{-kM}  \|F_{m}^{l}\|_{L^2}\|F_{m^{\prime}}^{l}\|_{L^2}.
\end{align}

For fixed $Q,Q^{\prime}: cQ \cap cQ^{\prime} \neq \emptyset$,  we may assume that $cQ \cap cQ^{\prime}$ is contained in  a $2^{k-j} \times 2^{k-j}$ cube centered at the origin, then
\begin{align}\label{Eq3.10++}
& \biggl\| \mathop{sup}_{t \in (0,1)}\biggl| \int_{\mathbb{R}^{2}} K(x,z,t) \Pi_{Q}(\frac{z}{2^{k}}) F_{m}^{l}(z)dz \biggl|  \biggl| \int_{\mathbb{R}^{2}} K(x,z^{\prime},t) \Pi_{Q^{\prime}}(\frac{z^{\prime}}{2^{k}}) F_{m^{\prime}}^{l}(z^{\prime})dz^{\prime} \biggl|  \biggl\|_{L^{1}(cQ \cap cQ^{\prime} )} \nonumber\\
&\le 2^{2k-2j}\biggl\| \mathop{sup}_{t \in (0,1)}\biggl| T_{2^{2k-j}}(\Pi_{Q}(\frac{\cdot}{2^{k}}) F_{m}^{l})(x,t) \biggl|  \biggl| T_{2^{2k-j}}(\Pi_{Q^{\prime}}(\frac{\cdot}{2^{k}}) F_{m^{\prime}}^{l})(x,t) \biggl|  \biggl\|_{L^{1}(B(0,1))}.
\end{align}
Notice that by uncertainty principle and the assumption that $2^{l} \ll 2^{(2k-j)/4}$, the Fourier transform of $\Pi_{Q}(\frac{\cdot}{2^{k}}) F_{m}^{l}$ and $\Pi_{Q^{\prime}}(\frac{\cdot}{2^{k}}) F_{m^{\prime}}^{l}$ are supported in $2^{-l}$ cubes with separation $\sim2^{-l}$, and $2^{l} \ll 2^{\frac{2k-j}{4}}$. Therefore, we can apply Theorem \ref{theorem1.5} to get
\begin{align}\label{Eq3.9+}
&\biggl\| \mathop{sup}_{t \in (0,1)}\biggl| T_{2^{2k-j}}(\Pi_{Q}(\frac{\cdot}{2^{k}}) F_{m}^{l})(x,t) \biggl|  \biggl| T_{2^{2k-j}}(\Pi_{Q^{\prime}}(\frac{\cdot}{2^{k}}) F_{m^{\prime}}^{l})(x,t) \biggl|  \biggl\|_{L^{1}(B(0,1))} \nonumber\\
&\le C2^{-\frac{l}{2}} 2^{(2k-j)(\frac{3}{4}+\epsilon)}  2^{-2(2k-j)}  \biggl\|\Pi_{Q}(\frac{\cdot}{2^{k}}) F_{m}^{l}\biggl\|_{L^{2}} \biggl\|\Pi_{Q^{\prime}}(\frac{\cdot}{2^{k}}) F_{m^{\prime}}^{l}\biggl\|_{L^{2}}.
\end{align}

Then inequalities (\ref{sum}) - (\ref{Eq3.9+})  imply that
\begin{align}
& \biggl\| \mathop{sup}_{t \in I}|e^{it\Delta}f_{m}^{l}| | e^{it\Delta}f_{m^{\prime}}^{l}| \biggl\|_{L^{1}(B(0,1))} \nonumber\\
&\leq C 2^{2k-2j}2^{-\frac{l}{2}} 2^{(2k-j)(\frac{3}{4}+\epsilon)}  2^{-2(2k-j)} \sum_{Q,Q^{\prime}: cQ \cap cQ^{\prime} \neq \emptyset} \biggl\|\Pi_{Q}(\frac{\cdot}{2^{k}}) F_{m}^{l}\biggl\|_{L^{2}} \biggl\|\Pi_{Q^{\prime}}(\frac{\cdot}{2^{k}}) F_{m^{\prime}}^{l}\biggl\|_{L^{2}} \nonumber\\
 &\quad +C  2^{-kM}  \|F_{m}^{l}\|_{L^2}\|F_{m^{\prime}}^{l}\|_{L^2} \nonumber\\
  & \le  C 2^{2k-2j}2^{-\frac{l}{2}} 2^{(2k-j)(\frac{3}{4}+\epsilon)}  2^{-2(2k-j)}\|F_{m}^{l}\|_{L^2}\|F_{m^{\prime}}^{l}\|_{L^2} \nonumber\\
 & =  C 2^{2k-2j}2^{-\frac{l}{2}} 2^{(2k-j)(\frac{3}{4}+\epsilon)}  2^{-2(2k-j)}2^{2k}\|f_{m}^{l}\|_{L^2}\|f_{m^{\prime}}^{l}\|_{L^2} \nonumber\\
  & \le C  2^{-\frac{l}{2}} 2^{(2k-j)\frac{3}{4}}2^{2k\epsilon} \|f_{m}^{l}\|_{L^2}\|f_{m^{\prime}}^{l}\|_{L^2}.
\end{align}
This completes the proof of Lemma \ref{main}.


\begin{flushleft}
\vspace{0.3cm}\textsc{Wenjuan Li\\School of Mathematics and Statistics\\Northwest Polytechnical University\\710129\\Xi'an, People's Republic of China}

\vspace{0.3cm}\textsc{Huiju Wang\\School of Mathematics Sciences\\University of Chinese Academy of Sciences\\100049\\Beijing, People's Republic of China}

\vspace{0.3cm}\textsc{Dunyan Yan\\School of Mathematics Sciences\\University of Chinese Academy of Sciences\\100049\\Beijing, People's Republic of China}

\end{flushleft}

\end{document}